# AN ALGORITHMIC AND A GEOMETRIC CHARACTERIZATION OF COARSENING AT RANDOM

By Richard D. Gill[1] and Peter D. Grünwald[2,3]

*Leiden University and CWI*

We show that the class of conditional distributions satisfying the coarsening at random (CAR) property for discrete data has a simple and robust algorithmic description based on randomized uniform multicovers: combinatorial objects generalizing the notion of partition of a set. However, the complexity of a given CAR mechanism can be large: the maximal "height" of the needed multicovers can be exponential in the number of points in the sample space. The results stem from a geometric interpretation of the set of CAR distributions as a convex polytope and a characterization of its extreme points. The hierarchy of CAR models defined in this way could be useful in parsimonious statistical modeling of CAR mechanisms, though the results also raise doubts in applied work as to the meaningfulness of the CAR assumption in its full generality.

**1. Introduction.** In statistical practice one is often presented with incomplete, or more generally, *coarse* data. To properly model such data, one needs to take into account the mechanism by which the data are coarsened. In practice the details of this coarsening mechanism are often unknown or computationally expensive to model. Therefore, it is of interest to determine conditions under which this mechanism can be safely ignored. The "coarsening at random" (CAR) assumption is the weakest condition giving this guarantee. It was identified by Heitjan and Rubin (1991). More recently,

Received June 2007; revised June 2007.
[1]Supported by the IST Programme of the European Community within the PASCAL Network of Excellence, IST-2002-506778. RDG's work was supported (during his previous affiliation) by the Department of Mathematics, Utrecht University, and, thanks to a visiting position, by the Thiele Centre, Århus University.
[2]CWI is the Dutch national research institute for mathematics and computer science.
[3]EURANDOM is funded by the Dutch science foundation, NWO, and Eindhoven University.

*AMS 2000 subject classifications.* Primary 62A01; secondary 62N01.
*Key words and phrases.* Coarsening at random (CAR), ignorability, uniform multicover, Fibonacci numbers.

Grünwald and Halpern (2003) and Jaeger (2005b) have stressed that the importance of CAR is not restricted to statistical applications: when updating a probability distribution based on new information in learning, artificial intelligence, or other scientific applications, it precisely characterizes when one can ignore the distinction between the fact that an event has been *observed*, and the fact that an event has *happened*, thereby considerably simplifying the update process.

Thus, both in statistical inference with coarsened data and for probability updating in learning algorithms, it is attractive to be able to make the CAR assumption. In order to be able to judge whether or not the assumption is warranted, it is important to fully understand its meaning. Here we approach this problem by giving two intimately related characterizations of the CAR assumption. First, we show that the set of all CAR mechanisms for a given finite sample space can be seen as a convex polytope. Each CAR mechanism is a mixture of CAR mechanisms which correspond to the vertices of the polytope. Our first main result, Theorem 1, characterizes these vertices. Our second result, which follows easily from the first, complements this geometric view with an algorithmic one. We show that a simple probabilistic algorithm can simulate any possible CAR mechanism, and only CAR mechanisms. Prompted by Gill, van der Laan and Robins (1997), earlier authors [Grünwald and Halpern (2003) and Jaeger (2005b)] have also searched for such constructions, calling them *procedural models* for CAR. Yet the procedural models proposed so far are not quite satisfactory, because in all cases,

1. The procedural model depends on parameters which have to be fine-tuned in order to guarantee the CAR property; or equivalently,
2. A small perturbation in the parameters can destroy the CAR property.

This "frailty" or lack of robustness is an indication that such procedures may not occur naturally. In fact Jaeger (2005b), Theorem 4.17, shows that the only CAR mechanisms which a robust procedure can generate must be of a special type known as "coarsening completely at random," CCAR.

Here we present a natural way to generate all CAR mechanisms, and only CAR mechanisms, *that does not require fine-tuning of parameters*. Our algorithm works for arbitrary finite sample spaces. It is based on a generalization of the notion of a partition of a set which we call a *uniform multicover*, or just multicover for short.

Superficially, its existence would have to contradict Jaeger's theorem mentioned above. But of course, a proven theorem does not allow any contradictions. The difference lies in the notion we use of robustness and of its negation, frailty. Our result can be seen as criticism of Jaeger's notion of robustness, even though this does at first sight seem appealing and natural. By parameterizing CAR distributions in a different manner, we obtain a representation in which CAR *can* be generated without parameter tuning. In a nut-shell: we consider a discrete uniform distribution to be a robust and natural object. Jaeger considers it to be an easily perturbed object.

We emphasize that the body of Jaeger's work remains highly relevant; this is just one of a number of important results he has obtained, and we, too, come to the conclusion that CAR mechanisms which are not CCAR will be very rare in practice. For instance, our final result, Theorem 3, shows that, although no fine-tuning is needed, the complexity (defined in terms of the "height" of multicovers) of the CAR mechanisms generated by our algorithm can grow exponentially in the size of the sample space.

The paper is organized as follows. In Section 2 we briefly introduce coarsening at random and other preliminaries. In Section 3 we give our geometric interpretation of CAR distributions (Theorem 1). In Section 4 we define uniform multicovers and use these to define our procedural CAR model. We show that it generates all and only CAR mechanisms (Theorem 2). In Section 5 we discuss our CAR model in detail. We show (Theorem 3) that it gives rise to an exponential lower bound on the height of the multicovers needed in Theorem 2. The proofs are given in the final section.

**2. Preliminaries.** Let $E$ be a finite nonempty set, containing $n$ elements. A coarsening mechanism is a probabilistic rule which replaces any point $x$ in $E$ with a subset $A$ of $E$ containing $x$. Thus a coarsening mechanism is specified by a collection of (conditional) probabilities $\pi_A^x$ such that for all $x$, $\sum_{A \ni x} \pi_A^x = 1$. Intuitively, $x$ is generated by some process which for simplicity we will refer to as "Nature." But rather than observing $x$ directly, the statistician observes a coarsening of $x$, that is, a set $A$ containing $x$. We call $x$ the *underlying outcome* and $A$ the corresponding *observation*. The coarsening mechanism determines the $A$ that is observed given $x$; $\pi_A^x$ is the probability of observing the set $A$ with $A \ni x$, given that Nature has generated $x$. We define the *support* of such a coarsening mechanism as the set of $A \subseteq E$ for which $\pi_A^x > 0$ for some $x \in E$.

A coarsening mechanism satisfies the CAR (coarsening at random) property if and only if for all $x, x' \in A$,

(2.1) $$\pi_A^x = \pi_A^{x'} = \pi_A, \qquad \text{say.}$$

Intuitively, this means that the probability of observing $A$ is the same for all $x$ that are contained in $A$: the coarsening is done "at random," independently of the underlying $x$. We note that (2.1) is the definition of CAR employed by Gill, van der Laan and Robins (1997). It is called "strong CAR" by Jaeger (2005a). The definition is explained in detail by Gill, van der Laan and Robins (1997) and Jaeger (2005a); motivation, practical relevance and applications of the CAR property are discussed extensively by Gill, van der Laan and Robins (1997) and Grünwald and Halpern (2003).

Definition (2.1) shows that a CAR mechanism is specified by a collection of probabilities $\pi_A$ indexed by the nonempty subsets $A$ of $E$ satisfying

$$\sum_{A \ni x} \pi_A = 1 \qquad \forall x \in E. \tag{2.2}$$

We can therefore represent a CAR mechanism by the vector $\boldsymbol{\pi} = (\pi_A : \varnothing \subset A \subseteq E)$, where we assume the subsets $A$ to be ordered in some standard manner. For a given finite set of CAR mechanisms $\boldsymbol{\pi}_1, \ldots, \boldsymbol{\pi}_p$, and any probability vector $\boldsymbol{\lambda} = (\lambda_1, \ldots, \lambda_p)$, we define their *mixture* $\boldsymbol{\pi}' = \lambda_1 \boldsymbol{\pi}_1 + \cdots + \lambda_p \boldsymbol{\pi}_p$. The following two observations are immediate:

1. For each partition of $E$, there is a unique CAR mechanism that has exactly that partition as its support (for each set $A$ in the partition, $\pi_A^x = \pi_A = 1$, for all $x \in A$).
2. Each finite mixture of CAR mechanisms again represents a CAR mechanism.

These two observations suggest a simple procedural CAR model: Fix some integer $p > 0$ and pick $p$ (arbitrary) partitions $\mathcal{E}_1, \ldots, \mathcal{E}_p$ of $E$. Each of these induces a unique corresponding CAR mechanism. Now fix an arbitrary distribution $\boldsymbol{\lambda} = \lambda_1, \ldots, \lambda_p$ on $\mathcal{E}_1, \ldots, \mathcal{E}_p$. The coarsened data are now generated by first, independently of the underlying $x$, selecting one of the $p$ partitions according to the distribution $\boldsymbol{\lambda}$. Then, within the chosen partition, the unique $A$ is generated which contains the underlying $x$. One can think of each partition as a "sensor" with the help of which the data are observed. The procedure amounts to selecting a sensor completely at random, independently of the underlying $x$ generated by Nature. This procedural CAR model is called the CARGEN procedure by Grünwald and Halpern (2003). The "parameters" of this procedure are the number of partitions $p$, the partitions $\mathcal{E}_1, \ldots, \mathcal{E}_p$ and the distribution $\boldsymbol{\lambda}$. Clearly, for every setting of the parameters, the resulting algorithm defines a CAR mechanism. One may be tempted to think that, by an appropriate setting of the parameters, *all* CAR mechanisms can be simulated by CARGEN, but the following example shows that this is not the case:

EXAMPLE 1 [Gill, van der Laan and Robins (1997)]. Let $E = \{1, 2, 3\}$, $A_{12} = \{1, 2\}$, $A_{23} = \{2, 3\}$ and $A_{31} = \{3, 1\}$. Consider the coarsening mechanism $\boldsymbol{\pi}^*$ defined by

$$\pi_{A_{12}}^{*1} = \pi_{A_{12}}^{*2} = \pi_{A_{23}}^{*2} = \pi_{A_{23}}^{*3} = \pi_{A_{31}}^{*3} = \pi_{A_{31}}^{*1} = \tfrac{1}{2}, \tag{2.3}$$

and $\pi_A^{*x} = 0$ for all other $x \in E, A \subseteq E$. By (2.1) it is immediately seen that this is a CAR mechanism. But because the support of the mechanism is not a union of partitions of $E$, it cannot be simulated by the CARGEN procedure.

The example shows that the CARGEN procedure is incomplete: there exist CAR mechanisms which cannot be represented by any parameter setting of CARGEN. The question is now whether there exist "natural" procedural CAR models which are complete. In previous work, two candidates for such models were proposed: Grünwald and Halpern's (2003) CARGEN$^*$ (an extension of CARGEN described above) and Jaeger's (2005b) *Propose-and-Test* model. Both of these suffer from the frailty property mentioned in the Introduction: rather than producing CAR mechanisms for all parameter settings, the parameters need to be fine-tuned. In previous work, one other procedural model has been proposed which, like CARGEN, produces CAR mechanisms for all settings of its parameters. However, as shown by Jaeger (2005b), this *randomized monotone coarsening* model [Gill, van der Laan and Robins (1997)] is in fact equivalent to CARGEN: both can simulate exactly the set of "coarsening completely at random" (CCAR) mechanisms. In fact, [Jaeger (2005b), Theorem 4.17] shows that any CAR mechanism that is not CCAR is, in a certain sense, nonrobust. For the details of Jaeger's definition of robustness we refer to Jaeger (2005b). Briefly, he supposes that a CAR mechanism involves an auxiliary randomization, and defines robustness in terms of robustness to changes in the distribution of the auxiliary variable.

Jaeger's result suggests that there exists no procedural CAR model that is both complete and does not require any parameter tuning. Yet in Section 4, we exhibit a simple extension of the CARGEN procedure which achieves exactly this, as long as we are able to sample from a uniform distribution. The procedural model will be based on a geometric interpretation of CAR which we present below.

**3. A geometric view of CAR.** We have already indicated that a finite mixture of CAR mechanisms $\pi$ is itself a CAR mechanism. Hence, for a given finite sample space $E$ the set of all CAR mechanisms defined with respect to $E$ forms a convex body in Euclidean space. In Theorem 1 we show that this body is a polytope with a finite number of extreme points, the vertices of the polytope. In order to characterize these extreme points, we first note that the support of a CAR mechanism is always a cover of $E$. With any cover of $E$ we associate its incidence matrix: the matrix $M$ with rows indexed by $x \in E$, columns indexed by $A$ in the support, and elements $\mathbb{1}_{\{x \in A\}}$. An incidence matrix of a cover is a matrix of 0's and 1's with at least one 1 in every row and column. We now use these incidence matrices to define extreme CAR mechanisms in an algebraic way. Theorem 1 below states that these CAR mechanisms are also extreme points in the geometric sense, justifying our terminology.

In the sequel, vectors are always column vectors, even if we lazily list the elements in a row. **0** and **1** denote vectors of 0's and 1's, respectively, whose length depends on the context.

Take the incidence matrix $M$ of an arbitrary cover $(A_1, \ldots, A_m)$ of $E$. If the equation $M\mathbf{z} = \mathbf{1}$ has a nonnegative solution, then this solution $\mathbf{z} = (z_1, \ldots, z_m)$ represents a CAR mechanism $\boldsymbol{\pi}$, where for any $A_j$ appearing in the cover, $z_j = \pi_{A_j}$, and for any $A$ not appearing in the cover, $\pi_A = 0$ [see also Grünwald and Halpern (2003), who explain this in detail]. We call $\boldsymbol{\pi}$ a CAR mechanism corresponding to $M$.

DEFINITION 1. We call $\boldsymbol{\pi}$ an *extreme CAR mechanism* if it corresponds to an incidence matrix $M$ of a cover $(A_1, \ldots, A_m)$ such that $M\mathbf{z} = \mathbf{1}$ has a unique, and strictly positive, solution.

By definition, a CAR mechanism is extreme if and only if it is the *only* CAR mechanism with the same support. It is easily checked that the mechanism $\boldsymbol{\pi}^*$ of Example 1 is an example of an extreme CAR mechanism: it is the only CAR mechanism with support $A_{12}, A_{23}, A_{31}$. The uniqueness also implies that the support of an extreme CAR mechanism cannot have more than $n$ elements (the size of $E$). It is clear that the number of extreme CAR mechanisms, for given $E$, is finite. We can find them all by enumerating and testing all covers of $E$ with $m \leq n$ elements.

THEOREM 1. *Every CAR mechanism is a mixture of extreme CAR mechanisms.*

In other words, all CAR mechanisms can be represented by randomly choosing, independently of $x$, one of a finite set of extreme CAR mechanisms. In the next section, we show that all such extreme mechanisms are of a simple and natural form. This will lead to Theorem 2, a direct corollary of Theorem 1, giving an algorithmic characterization of CAR.

**4. An algorithmic view of CAR.** Our procedure is based on the notion of a *uniform multicover*, which we now define. A *k-multicover* of $E$, or just $k$-cover for short, is a collection of nonempty subsets of $E$, allowing multiplicities, such that for each $x \in E$, precisely $k$ of the sets (some of which may be the same) contain $x$. Thus a 1-cover is an ordinary partition of $E$. By a *uniform multicover* we mean a $k$-cover for some $k \geq 1$. The *height* of a uniform multicover is its value of $k$. The *support* of a multicover is the set of subsets of $E$ in the multicover.

A $k$-cover is specified by its support and by the multiplicity of each set in its support. Thus, to each nonempty subset $A$ of $E$ there corresponds a nonnegative integer $n_A$ such that $n_A = 0$ if $A$ is absent from the $k$-cover; otherwise $n_A > 0$ is the multiplicity of $A$ in the $k$-cover. The $n_A$ have to satisfy

$$\sum_{A \ni x} n_A = k \qquad \forall x \in E. \tag{4.1}$$

For a given $k$-cover we can now define a CAR mechanism by setting

(4.2) $$\pi_A = n_A/k \qquad \forall A \subseteq E.$$

The algorithmic interpretation is as follows: Nature generates some $x \in E$. The coarsening mechanism investigates which $A$ in the uniform multicover contain $x$. There are exactly $k$ such $A$, including multiplicities, whatever $x$. We choose one of these uniformly at random, that is, each $A$ with $x \in A$ is chosen with probability $1/k$.

Conversely, any CAR mechanism for which all the CAR probabilities $\pi_A$ are rational numbers is generated by a $k$-cover with $k$ equal to the lowest common multiple of the denominators of the $\pi_A$. We call CAR mechanisms obtained in this way *rational*. *The rational CAR mechanisms are precisely the CAR mechanisms generated by a uniform multicover.* Note that if $k$ and all $n_A$ share a common factor, we can divide by this factor without changing the $\pi_A$. We consider such multicovers as equivalent and take the multicover with the smallest $k$ as representative of the class. In this way, each rational CAR mechanism corresponds to exactly one uniform multicover, and vice versa. We can make the connection to Theorem 1 by noting that

FACT 1. *Every extreme CAR mechanism is rational. Thus, it is generated by a uniform multicover.*

This follows directly from the fact that the matrix $M$ in Definition 1 is a 0/1-matrix and the solution of $M\mathbf{z} = \mathbf{1}$ is unique.

As stated above, for each rational CAR mechanism there is a unique uniform multicover which generates it. We can thus define an "extreme multicover" as a uniform multicover that generates an extreme CAR mechanism. Using Theorem 1, it is easily shown that extreme multicovers are just those uniform multicovers that do not contain a subset that is also a uniform multicover (we omit the details of the reasoning).

We may now define a procedural CAR model by first fixing a finite number $p$ of arbitrary uniform multicovers $\mathcal{C}_1, \ldots, \mathcal{C}_p$. We then fix an arbitrary distribution $\boldsymbol{\lambda} = (\lambda_1, \ldots, \lambda_p)$ on $\mathcal{C}_1, \ldots, \mathcal{C}_p$. The coarsened data are now generated by first, independently of the underlying $x$, selecting one of the $p$ uniform multicovers according to the distribution $\boldsymbol{\lambda}$. Suppose we have chosen multicover $\mathcal{C}_j$ with height $k_j$. Then among the $k_j$ sets in $\mathcal{C}_j$ which contain $x$, we choose one uniformly at random, with probability $1/k_j$. This procedural CAR model is a simple extension of CARGEN (Section 2), where the role of partitions is taken over by the more general uniform multicovers. Like CARGEN, it simulates a CAR mechanism for all parameter settings; no fine-tuning is needed. Theorem 2(ii) below (a corollary of Theorem 1) states that by appropriately setting the parameters, we can simulate *all* CAR mechanisms. Before presenting the theorem, we continue our example.

EXAMPLE 2 (Example 1 continued). The collection $\mathcal{C} = \{A_{12}, A_{23}, A_{31}\}$ is a uniform multicover of $E$ with height 2. Consider a simple instantiation of the procedural CAR model we described above, with just one multicover $\mathcal{C} = \mathcal{C}_1$, so that $\boldsymbol{\lambda} = (1)$. For each $x$ chosen by Nature, there will be exactly two elements of $\mathcal{C}$ which contain $x$. We select between these with probability $1/2$. It is immediately clear that this algorithm simulates the CAR mechanism $\boldsymbol{\pi}^*$ described in Example 1. An implementation of this mechanism requires a fair coin toss. If the coin is biased the CAR property can be lost. Relatedly, the mechanism is not robust in Jaeger's sense.

THEOREM 2. (i) *Every CAR mechanism can be arbitrarily well approximated by a rational CAR mechanism, that is, for all CAR mechanisms $\boldsymbol{\pi}$, all $\varepsilon > 0$, there exists a rational CAR mechanism $\boldsymbol{\pi}'$ such that $\|\boldsymbol{\pi} - \boldsymbol{\pi}'\| < \varepsilon$.*

(ii) *Every CAR mechanism is exactly equal to a finite mixture of extreme (and hence rational) CAR mechanisms.*

We extensively discuss this theorem in the next section.

**5. Discussion.** Theorem 2 shows that there is an easy probabilistic algorithm which approximates each CAR mechanism arbitrarily well, and that a randomized version of the algorithm reproduces each one exactly. Since the rational numbers form a dense subset of the reals, Theorem 2(i) is, in a sense, trivial. The real innovation is part (ii), which shows that each CAR distribution can be represented *exactly* as a mixture of a finite set of candidate rational mechanisms.

No fine-tuning of parameters is required to ensure the CAR properties so the algorithms do have a robustness property. We just need to be able to choose uniformly at random from a finite set. Of course, if one perturbs the uniform distribution over the $k$ sets containing a point $x$, one will in general destroy the CAR property—this is the reason that our result does not contradict Jaeger's (2005b), Theorem 4.17. For this reason, some readers may not want to call the procedure "robust." However, the (weaker) claim that the algorithm requires no parameter tuning seems indisputable: we can hardly think of implementing a uniform distribution as "parameter tuning." Unlike the parameters in earlier complete procedural CAR models, which could vary from situation to situation and were hard to determine, the uniform distribution is universal and easy to determine. If the device we use to generate a uniform distribution does not work perfectly, our procedural model will slightly violate CAR, hence one might perhaps say it is "nonrobust"; but devices used to generate a uniform distribution (coins, dice) exist, and usually do not arise as fine-tuned versions of devices that can generate a whole range of distributions; hence one cannot say that our model requires "fine-tuning."

The reason that earlier complete procedural CAR models did require parameter tuning was that their parameters had to satisfy complicated constraints [see, e.g., Example 4.7 in Jaeger (2005b)]. As remarked by M. Jaeger, we do pay a price for avoiding these parameter constraints: we now have complicated constraints (4.1) on multiplicities of sets appearing in multicovers. Such constraints are arguably more natural than constraints on continuous-valued parameters, at least as long as the multicovers involved are not too complex. Unfortunately, in order to span all CAR mechanisms, we sometimes need highly complex multicovers, as we show below. This limits the importance of our procedural model, as we discuss further below.

We can measure the complexity of multicovers in terms of their height. Since the row rank of $M$ equals its (full) column rank, $m$, we can delete rows obtaining an $m \times m$ nonsingular matrix $M_0$. Deleting the corresponding rows from $\mathbf{1}$ also, we obtain $\mathbf{z} = M_0^{-1}\mathbf{1}$. It follows by the standard expression of matrix inverse in terms of determinants that the value of $k$ appearing in (4.2) is bounded by $m!$. Hence, the height of the extreme multicovers that can be defined on a sample space of size $|E| = n$ is upper bounded by $n!$. But is this too pessimistic? Unfortunately not, or at least, not significantly: our next and last theorem gives an exponential *lower* bound on the *maximal* height of an extreme multicover. It turns out that this grows at least as fast as the celebrated Fibonacci numbers, defined as $F_1 = 1, F_2 = 1$, and for $j \geq 3$, $F_j = F_{j-1} + F_{j-2}$.

Theorem 3 below considers $n \times n$ matrices $S_n$ inductively defined as follows: $S_1 = (1)$. For odd $n$, $S_{n+1}$ is constructed from $S_n$ by setting

$$S_{n+1} = \begin{pmatrix} 1 & \mathbf{0}^\top \\ \mathbf{0} & S_n \end{pmatrix}.$$

For even $n$, $S_{n+1}$ is constructed from $S_n$ by setting

$$S_{n+1} = \begin{pmatrix} 0 & \mathbf{1}^\top \\ \mathbf{1} & S_n \end{pmatrix}.$$

This is easier than it seems: the pattern should be obvious from the example $n = 9$, shown in Figure 1.

THEOREM 3. *For odd $n > 0$, the equation $S_n \mathbf{z} = \mathbf{1}$ has the unique solution*

$$\mathbf{z} = \left(\frac{F_{n-1}}{F_n}, \frac{F_{n-2}}{F_n}, \ldots, \frac{F_2}{F_n}, \frac{F_1}{F_n}, \frac{1}{F_n}\right),$$

*so that $S_n$ represents an extreme point for sample spaces with size $|E| = n$, with height $k = F_n$.*

The theorem implies that the maximal height of an extreme multicover grows exponentially fast with $n$; also, the maximal needed multiplicity of a set in an extreme multicover grows exponentially fast with $n$. We interpret this result as follows.

Uniform multicovers are important in two ways:

1. They lead to an attractive algorithmic characterization of CAR that requires no fine-tuning of parameters (Theorem 2).
2. They induce a hierarchy of CAR models that could be of use in statistical applications. We elaborate on this below.

Yet apart from these applications, the importance of uniform multicovers in understanding CAR is limited—the maximal needed height of the multicover grows exponentially fast with $n$, so though the idea of the algorithm is simple, its detailed specification can be complex. Thus, we can say neither that our characterization provides a truly simple description of every CAR mechanism, nor that our multicover CAR mechanisms always correspond to some "natural" process. While it seems reasonable to suppose that low-height multicovers may be good models for some processes occurring in nature, the same cannot be said for exponentially high multicovers, and our Theorem 3 does show that we need to take these into account.

Jaeger's (2005b) robustness Theorem 4.17 suggests that the CAR mechanisms occurring in nature are those generated by randomized 1-covers. Our characterization nuances this somewhat, suggesting that in some situations $k$-covers for small $k > 1$ may also be reasonable models. Indeed, the hierarchy of CAR mechanisms induced by our algorithm suggests a statistical estimation procedure for parsimoniously estimating CAR mechanisms and their parameters. Such a procedure would penalize the fit of a proposed CAR mechanism to the data. The penalization would be some function of the number of extreme multicovers needed to express the mechanism, and the height of each of these. Alternatively one could use just one multicover,

$$\begin{pmatrix} 0 & 1 & 1 & 1 & 1 & 1 & 1 & 1 & 1 \\ 1 & 1 & 0 & 0 & 0 & 0 & 0 & 0 & 0 \\ 1 & 0 & 0 & 1 & 1 & 1 & 1 & 1 & 1 \\ 1 & 0 & 1 & 1 & 0 & 0 & 0 & 0 & 0 \\ 1 & 0 & 1 & 0 & 0 & 1 & 1 & 1 & 1 \\ 1 & 0 & 1 & 0 & 1 & 1 & 0 & 0 & 0 \\ 1 & 0 & 1 & 0 & 1 & 0 & 0 & 1 & 1 \\ 1 & 0 & 1 & 0 & 1 & 0 & 1 & 1 & 0 \\ 1 & 0 & 1 & 0 & 1 & 0 & 1 & 0 & 1 \end{pmatrix}$$

Fig. 1. *The matrix $S_9$, an example of the matrices $S_n$ figuring in Theorem 3.*

not necessarily extreme, and penalize its height. This could be done either explicitly, by adding a regularization term to the likelihood, or implicitly, by the use of suitable Bayesian priors.

Such procedures could be useful in practice if one seriously believed that the data are CAR but quite possibly, not CCAR. One could hope in this way to combine the advantages of asymptotic validity and even go for asymptotic efficiency, with good small sample behavior. However, our results can also be read in a different way. Though we found an appealing way to model CAR, it remains the fact that there do not seem to be so many good reasons in practice, in general, to assume CAR but not CCAR. Therefore, if one is prepared to assume CAR, one is likely to be also prepared to assume CCAR. Though the distinction concerns a "nuisance" part of the model, and indeed, in likelihood approaches is invisible by the likelihood factorization implied by CAR, one can capitalize on the extra knowledge, for instance, in order to obtain better small sample properties of estimators, at the cost of loss of asymptotic efficiency.

A final view is that the extra generality obtained by relaxing CCAR to CAR is illusory. If one does not believe in CAR, one has no option but to start modeling and estimating the coarsening mechanism. Jaeger (2006a, 2006b) has made some proposals in this direction which seem promising. Another possibility, so far not explored, is to use the notion of *relative* rather than *absolute* CAR introduced by Gill, van der Laan and Robins (1997). The point of CAR is that, in likelihood inference, one can analyze coarsened data *as if* the coarsening mechanism had been fixed in advance as any particular CAR mechanism, and specifically therefore, as if coarsening by an independently fixed-in-advance partitioning of the sample space. Relative CAR means CAR relative to some other specific (non-CAR) coarsening mechanism: the likelihood factors; the interesting part is the same *as if* the data had been coarsened by the reference coarsening model; the nuisance part can be used for inference concerning *which* coarsening mechanism has generated the data, out of the mechanisms in the family implied by the reference mechanism. It would be interesting to explore this possibility in more detail.

## 6. Proofs.

6.1. *Proof of Theorem 1.* We show below that the set of all CAR mechanisms forms a convex polytope and characterize the extreme points in terms of linear algebra, corresponding to Definition 1.

A CAR mechanism is a collection of numbers $\pi_A$ indexed by the nonempty subsets $A$ of a finite set $E$. They must satisfy two sets of constraints: the inequalities $\pi_A \geq 0$ for each $A$, and the equalities $\sum_{A \ni x} \pi_A = 1$ for each $x$, both of which are obviously linear. Together the constraints imply that $\pi_A \leq 1$ for all $A$. Collecting the $\pi_A$ into a vector $\boldsymbol{\pi}$, we see that the set of all $\boldsymbol{\pi}$ is a convex, compact polytope since it is bounded and is the intersection of a finite number of closed half-spaces (one for each inequality constraint) and hyperplanes (one for each equality constraint). Hence each $\boldsymbol{\pi}$ is a convex combination of the extreme points of the polytope, of which there are a finite number in total.

The polytope lives in the affine subspace of all vectors $\boldsymbol{\pi}$ satisfying the equality constraints $\sum_{A \ni x} \pi_A = 1$ for each $x$. Since $\boldsymbol{\pi}$ has $2^{|E|} - 1$ components (the number of nonempty subsets of $E$) and there are $|E|$ constraints, it follows that the dimension of this affine subspace is $2^{|E|} - 1 - |E|$. The polytope is just the intersection of that affine subspace with the positive orthant. Within the affine subspace, each face of the polytope corresponds to one of the hyperplanes $\pi_A = 0$. Each vertex of the polytope is the unique meeting point of a number of faces; one for each $A$ such that $\pi_A = 0$. Thus to each vertex is associated a collection of subsets $A$ such that if we set the corresponding $\pi_A$ equal to 0 in the equations $\sum_{A \ni x} \pi_A = 1$ for all $x$, there is a unique and strictly positive solution in the remaining $\pi_A$. Conversely, any such collection of $A$ defines a vertex.

The subsets $A$ *not* in the collection define the support of the extreme CAR mechanism $\boldsymbol{\pi}$ under consideration. Let $M$ be its incidence matrix: the matrix of zeros and ones with rows indexed by elements $x \in E$, columns indexed by $A$ in the support, and with entries $\mathbb{1}_{\{x \in A\}}$. Write $\boldsymbol{\pi}_0$ for the vector of $\pi_A$ for $A$ in the support. In matrix form, the equations which must have a unique and positive solution $\mathbf{z} = \boldsymbol{\pi}_0$ can be written

$$M\mathbf{z} = \mathbf{1}, \tag{6.1}$$

and we have proved that there is a one-to-one correspondence between vertices of the polytope and incidence matrices $M$ of covers of $E$ such that this equation has a unique and positive solution. As we argued in Section 4, if the solution is unique it has to be rational.

Combining these facts, extreme points of the polytope of CAR mechanisms correspond to covers of $E$ whose incidence matrix $M$ is such that $M\mathbf{z} = \mathbf{1}$ has a unique solution, and the solution is strictly positive.

REMARK. A condition equivalent to $M\mathbf{z} = \mathbf{1}$ having a unique positive solution [Farkas's lemma in the theory of linear programming Schrijver (1986), Chapter 7] is that $M$ has full column rank, and, if $\mathbf{y}$ is such that (a) $\mathbf{y}^\top M \geq \mathbf{0}$, then (b) $\mathbf{y}^\top \mathbf{1} \geq 0$, with equality in (b) implying equality in (a). By arguments from integer programming [see again Schrijver (1986)] one may restrict here to vectors $\mathbf{y}$ of integers. Jaeger (2005b) gives a version of this condition for the existence of a CAR mechanism with given support— he does not demand full rank since he does not ask for uniqueness. Though more combinatorial in nature, this version of the condition for extremality does not seem to be much more useful, except perhaps for helping one to show that certain covers do *not* lead to solutions.

6.2. *Proof of Theorem 2.* Theorem 2 is, in fact, a direct corollary of Theorem 1. Namely, each extreme point is rational and therefore corresponds to a uniform multicover. Every point in a polytope can be written as a mixture of its extreme points. This gives us item (ii). Item (i) follows by considering the rational convex combinations of the extremes, which lie dense in all convex combinations.

6.3. *Proof of Theorem 3.* We prove the theorem by induction on $n$. For $n = 1$, the result trivially holds. Now suppose the result holds for $S_{n-1}$, for some even $n > 1$. Thus, $S_{n-1}\mathbf{q} = \mathbf{1}$ has a unique solution

$$(6.2) \quad \mathbf{q} = (q_1, \ldots, q_{n-1}) = \left(\frac{F_{n-2}}{F_{n-1}}, \frac{F_{n-3}}{F_{n-1}}, \ldots, \frac{F_2}{F_{n-1}}, \frac{F_1}{F_{n-1}}, \frac{1}{F_{n-1}}\right).$$

We prove the theorem by showing that this implies that

$$(6.3) \quad S_{n+1}\mathbf{r} = \mathbf{1}$$

has the unique solution

$$(6.4) \quad \mathbf{r} = (r_1, \ldots, r_{n+1}) = \left(\frac{F_n}{F_{n+1}}, \frac{F_{n-1}}{F_{n+1}}, \ldots, \frac{F_2}{F_{n+1}}, \frac{F_1}{F_{n+1}}, \frac{1}{F_{n+1}}\right).$$

To prove (6.4), note first that to each row of (6.3) corresponds a linear equation. Writing the equations corresponding to the first two rows explicitly and the equations corresponding to rows 3 to $n+1$ in matrix form, and reordering terms, we see that (6.3) is equivalent to

$$(6.5) \quad r_2 = 1 - \sum_{i=3}^{n+1} r_i,$$

$$(6.6) \quad r_2 = 1 - r_1,$$

$$(6.7) \quad S_{n-1}(r_3, \ldots, r_{n+1})^T = 1 - r_1,$$

where, by our inductive assumption, the last equality implies

$$(6.8) \quad (r_3, \ldots, r_{n+1}) = (1 - r_1)(q_1, \ldots, q_i),$$

and in particular

$$(6.9) \quad \sum_{i=3}^{n+1} r_i = (1 - r_1) \sum_{i=1}^{n-1} q_i.$$

Combining (6.6) with (6.5), we get $r_1 = \sum_{i=3}^{n-1} r_i$. Plugging this into (6.9) gives

$$\frac{r_1}{1 - r_1} = \sum_{i=1}^{n-1} q_i \tag{6.10}$$

where $q_i$ are given by (6.2). We claim this has the unique solution $r_1 = F_n/F_{n+1}$. To see this, note the following basic fact which follows immediately from repeatedly substituting the definition $F_n = F_{n-1} + F_{n-2}$ on the left in (6.11):

FACT 2. *For odd $n > 0$,*

$$F_n = \sum_{i=1}^{n-2} F_j + 1. \tag{6.11}$$

The fact implies that the right-hand side of (6.10) is equal to $F_n/F_{n-2}$. Plugging in our proposed solution $r_1 = F_n/F_{n+1}$, the left-hand side of (6.10) becomes $F_n/(F_{n+1} - F_n) = F_n/F_{n-2}$, so that (6.10) holds. This shows that $r_1$ is indeed given by $F_n/F_{n+1}$. By (6.6) it now follows that $r_2 = F_{n-1}/F_{n+1}$, and, by (6.8), that for $j \in \{3, \ldots, n+1\}$, $r_j = q_{j-2}/r_2 = s_j/F_{n+1}$, where $(s_1, s_2, \ldots, s_{n-2}, s_{n-1}) = (F_{n-2}, F_{n-3}, \ldots, F_1, 1)$. This shows that (6.4) is the unique solution of $S_{n+1}\mathbf{r} = \mathbf{1}$, and thus completes the induction step. The theorem is proved.

**Acknowledgments.** We would like to thank Sasha Gnedin and Lex Schrijver for stimulating conversations. Lex Schrijver made some essential contributions to the proof of Theorem 3.

Mathematical Institute  
Leiden University  
Niels Bohrweg 1  
2333 CA Leiden  
Netherlands  
E-mail: [gill@math.leidenuniv.nl](gill@math.leidenuniv.nl)  
URL: [http://www.math.leidenuniv.nl/~gill](http://www.math.leidenuniv.nl/~gill)

CWI  
Kruislaan 413  
1096 SJ Amsterdam  
and  
Eurandom  
Eindhoven  
Netherlands  
E-mail: [peter.grunwald@cwi.nl](peter.grunwald@cwi.nl)  
URL: [http://www.grunwald.nl](http://www.grunwald.nl)